\documentclass[reqno,11pt]{amsart}
\usepackage{graphicx}
\usepackage{amsmath,amssymb}
\usepackage{eepic}
\usepackage{epic}
\usepackage{amsmath,amssymb}
=-2.5 cm \voffset =-2.0 cm \columnsep=0.8cm

\def\intl {\int\limits}

\def\intl{\int\limits}
\def\ov{\overline}

\begin{document}
\twocolumn[ \begin{center}{\bf HILBERT'S 16th PROBLEM FOR QUADRATIC
SYSTEMS.

NEW METHODS
 BASED ON A TRANSFORMATION TO THE LIENARD EQUATION}

\bigskip
  {G. LEONOV }

{St. Petersburg State University, Universitetsky prospekt, 28,}

{Peterhoff,  St. Petersburg, 198904, Russia\\
email: leonov@math.spbu.ru }
\end{center}
\bigskip

{\begin{quote}\textwidth=14.5cm

Fractionally-quadratic transformations which reduce any
two-dimensional quadratic system to the special Lienard equation are
introduced. Existence criteria of cycles are obtained.
\medskip

\end{quote}}
\bigskip ]

\bigskip

{\bf 1. Introduction}
\bigskip

Within last  century the study of cycles of two-dimensional
quadratic systems was stimulated by   16th Hilbert's problem and its
different variants [Hilbert, 1902; Lloyd, 1988; Blows \& Perko,
1994; Ilyashenko, 2002; Chen \& Wang, 1979; Shi, 1980].

The development of asymptotical methods for investigation of cycles
of quadratic systems    is facilitated by the fractionally-quadratic
transformations, which reduce any two-dimensional quadratic system
to the Lienard equation of special type with discontinuous nonlinear
functions [Leonov, 1997; Leonov, 1998; Leonov, 2006]. In the present
paper these new facilities are demonstrated.

The elements of    direct Lyapunov method, the     small
perturbations  methods, and the bifurcation theory together permits
us to obtain the existence criteria of small cycles in the
neighborhoods of stationary points of Lienard equation.

 This criteria applied to quadratic systems.

The question arises as to whether in two--dimen\-sional quadratic
systems a simultaneous bifurcation of cycles from two weak stable
and weak unstable equilibria is possible by variation of one scalar
parameter. Here we provide the positive answer to this question.

For the Lienard equation,  the classical theorems on the existence
of cycles is well known [Lefschetz, 1957; Cesari, 1959]. Here there
occurs quite naturally the possibility of generalization  of these
results and of applying them to the problem on the existence of
cycles in quadratic systems. On this way in the space of parameters
of quadratic systems, the sets of positive Lebesgue measure, for
which the cycles exist, are discovered.

\bigskip
\begin{center}

{\bf 2. The transformation of quadratic system to the Lienard
equation }
\end{center}
\bigskip

Let us consider the quadratic system
$$\begin{aligned}
&\dot x=a_1x^2+b_1xy+c_1y^2+\alpha_1x+\beta_1y\\
&\dot y=a_2x^2+b_2xy+c_2y^2+\alpha_2x+\beta_2y\end{aligned}\eqno(1)
$$
where $a_i,b_i,c_i,\alpha_i,\beta_i$ are real numbers.

Here we shall    follow the paper [Leonov, 1997; Leonov, 1998;
Leonov, 2006].

{\bf Proposition 1.} {\it Without loss of generality, we can assume that
$c_1=0$. }

\medskip
\def\k{\kappa}
\def\al{\alpha}
{\bf Proof.} Suppose, for definiteness, that $a_2\neq0$
(otherwise, having performed the change $x\to y$, $\,y\to x$,
we obtain at once $c_1=0$). Introduce the linear transformation
$x_1=x+\nu y$, $\ y_1=y$. To prove Proposition 1,
it is sufficient to show that for the certain numbers $\rho$, $\k$, $\nu$
the following identity
$$
\begin{aligned}
(x&+\nu y)^{\bullet}=(a_1+\nu a_2)x^2+(b_1+\nu b_2)xy+\\
&+ (c_1+\nu c_2)y^2+(\al_1+\nu\al_2)x+(\beta_1+\nu\beta_2)y=\\
&=\rho(x+\nu y)y+\k(x+\nu
y)^2+(\al_1+\nu\al_2)x+\\
&+(\beta_1+\nu\beta_2)y
\end{aligned}
$$
is valid.
This identity is equivalent to the following system of equations
$$
\begin{array}{l}
\k=a_1+\nu a_2,\cr \k\nu^2+\rho\nu=c_1+\nu c_2,\cr
\rho+2\k\nu=b_1+\nu b_2.
\end{array} \eqno(2)
$$
The above relations are satisfied if
$$
(a_1+\nu a_2)\nu^2-\nu(b_1+\nu b_2)+(c_1+\nu c_2)=0.
$$
Since $a_2\neq0$, this  equation with respect to $\nu$ always has a
real root. Thus, the system of equations (2)    always has a real
solution. \hfill$\Box$

Further we assume that $c_1=0$.

\medskip

{\bf Proposition 2.} {\it Let be $b_1\neq0$. The straight line
$\beta_1+b_1x=0$ on the plane $\{x,y\}$ is invariant or
transversal for trajectories of system (1). }

\smallskip

{\bf Proof.} The last statement follows from the relation
$$\begin{aligned}
&(\beta_1+b_1x)^{\bullet}{=}b_1\big[(b_1x+\beta_1)y+a_1x^2+\al_1x\big]
{=}\\
&= \left[a_1
\left(\frac{\beta_1}{b_1}\right)^2-\al_1\left(\frac{\beta_1}
{b_1}\right)\right]b_1\end{aligned}
$$
for $x=-\beta_1/b_1$.
This implies that if
$$
a_1\left(\frac{\beta_1}{b_1}\right)^2-\al_1\left(
\frac{\beta_1}{b_1}\right)=0,
$$
then the straight line $\beta_1+b_1x=0$ is invariant and if
$$
a_1\left(\frac{\beta_1}{b_1}\right)^2-\al_1\left(
\frac{\beta_1}{b_1}\right)\neq0,
$$
then the straight line $\beta_1+b_1x=0$ is transversal. \hfill$\Box$

        Excluding the trivial case     that the right-hand
side of the first equation of system (1) is independent of $y$, we
assume that
$$
|b_1|+|\beta_1|\neq0. \eqno(3)
$$

Then by Proposition 2 we obtain that the trajectories of system (1)
are also the trajectories of the system
$$
\begin{aligned}
&\dot x=y+\frac{a_1x^2+\al_1x}{\beta_1+b_1x},\\
 &\dot
y=\frac{a_2x^2+b_2xy+c_2y^2+\al_2x+\beta_2y} {\beta_1+b_1x}.
\end{aligned} \eqno(4)
$$

Consider the following transformation
$$
\begin{array}{l}
\bar y=y+\dfrac{a_1x^2+\al_1x}{\beta_1+b_1x},\cr \bar x=x.
\end{array}
$$
Using these new phase variables (here the bars over the variables
$\bar x\to x$, $\bar y\to y$ are omitted), system (4) can be written as
$$
\dot x=y,\qquad \dot y=-Q(x)y^2-R(x)y-P(x), \eqno(5)
$$
where
$$
Q(x)=\frac{-c_2}{\beta_1+b_1x},
$$

$$
\begin{aligned}
R(x)=-&\dfrac{(b_1b_2-2a_1c_2+a_1b_1)x^2+(b_2\beta_1+b_1\beta_2
-}{(\beta_1+b_1x)^2}\cr\cr
&\dfrac{-2\al_1c_2+2a_1\beta_1)x+\al_1\beta_1+\beta_1\beta_2}
{(\beta_1+b_1x)^2},
\end{aligned}
$$
$$\begin{aligned}
P(x)&=-\left(\frac{a_2x^2+\al_2x}{\beta_1+b_1x}-
\frac{(b_2x+\beta_2)(a_1x^2+\al_1x)}{(\beta_1+b_1x)^2}{+}\right.\\
&\left.+\frac{c_2(a_1x^2{+}\al_1x)^2}{(\beta_1{+}b_1x)^3}\right).\end{aligned}
$$

Proposition 2 and condition (3) imply that the trajectories of system (5)
are also trajectories of the system
$$
\dot x=ye^{p(x)},\quad\dot y=\big[-Q(x)y^2-R(x)y-P(x)\big]e^{p(x)},
$$
where $p(x)$ is a certain integral of the function $Q(x)$.

From this system, using the change $\bar x=x$, $\bar y=ye^{p(x)}$,
we obtain
$$
\dot x=y,\qquad \dot y=-f(x)y-g(x). \eqno(6)
$$
Here the bars over the variables $x$ and $y$ are also omitted.

So, for $b_1\ne0$, by means of the above nondegenerate changes
we can reduce system (1) to the Lienard equation (6)
with the functions
$$
f(x)=R(x)e^{p(x)}=R(x)|\beta_1+b_1x|^q,
$$
$$
g(x)=P(x)e^{2p(x)}=P(x)|\beta_1+b_1x|^{2q}
$$
Here for $b_1\ne0$,\quad$q=-\frac{c_2}{b_1}$.

{\bf Proposition 3.} {\it Let be $b_1\ne 0,\,\beta_1\ne
0,\,\,\alpha_1\ne 0$. Without loss of generality, we can assume that
$b_1=\alpha_1=\beta_1=1$.}

{\bf Proof.} Using the change
$$
x=\frac{\beta_1}{b_1}\ov{x},\,\,y=\frac{\alpha_1}{b_1}\ov{y},\,\,t=\frac{\ov{t}}{\alpha_1},
$$
we obtain
$$
\dot{\ov{x}}=\frac{\alpha_1\beta_1
a_1}{b_1}{\ov{x}}^2+\ov{x}\ov{y}+\ov{x}+\ov{y}.\eqno\Box
$$

Let us consider system (6) where
$$\begin{aligned}
&f(x)=(Ax+B)x|x+1|^{q-2},\\
&g(x)=\frac{(C_1x^3+C_2x^2+C_3x+C_4)x}{(x+1)^3}|x+1|^{2q}.
\end{aligned}
$$

Here $A,B,C_j\,\,(j=1,\ldots,4)$, $q$ are real numbers.

{\bf Proposition 4.} {\it For numbers $A,B,C_j,q$ of system (6)
exist corresponding numbers $a_1,
b_1=1,\,\alpha_1=1,\,\beta_1=1,\,a_2, b_2, c_2=-q,\,\alpha_2,
\beta_2=-1$ of system (1) iff
$$\begin{aligned}
&\frac{(B-A)}{(2q-1)^2}\left((1-q)B+(3q-2)A\right)=2C_2-3C_1-C_3,\\
&\frac{(B-A)}{(2q-1)^2}\left(B+2(q-1)A\right)=C_2-2C_1-C_4.\end{aligned}
$$
Here
$$\begin{aligned}
&a_1=1+\frac{B-A}{2q-1},\\
&a_2=-(q+1)a_1^2-Aa_1-C_1,\\
&b_2=-A-a_1(2q+1),\\
&\alpha_2=a_1^2-2a_1+A(a_1-1)+(2C_1-C_2).
\end{aligned}
$$}

\bigskip
\begin{center}
{\bf 3. Existence criterion for small cycles of the Lienard equation
in the neighborhood of equilibrium}

\end{center}

\bigskip

 Consider the equation
$$
\ddot z+z=u(t).\eqno(7)
$$

For this equation  it is well known the following simple result.

{\bf Lemma 1} [Arnol'd, 1976]. {\it The solution of equation (7)
with initial data $z(0), \dot z(0)$ is given by formula}
$$\begin{aligned}
z(t)&=\left[z(0)-\intl_0^tu(\tau)\sin\tau d\tau\right]\cos t+\\
&+\left[\dot z(0)+\intl_0^tu(\tau)\cos\tau d\tau\right]\sin t.
\end{aligned}\eqno(8)
$$

For $\dot z(0)=0$, taking into account formula (8),
we have the following relation
$$\begin{aligned}
\dot z(t)&=-z(0)\sin t+\intl_0^tu(\tau)\sin\tau d\tau\sin t+\\
&+\intl_0^tu(\tau)\cos\tau d\tau\cos t. \end{aligned}\eqno(9)
$$

Consider the equation
$$
\ddot x+F(x,\varepsilon)\dot x+G(x,\varepsilon)=0
$$
or the equivalent  system
$$\begin{aligned}
&\dot x=y\\
&\dot y=-F(x,\varepsilon)y-G(x,\varepsilon).
\end{aligned}\eqno(10)
$$
Here $\varepsilon$ is a positive number, $F(x,0)=f(x)$,
$G(x,0)=g(x)$, $f(x_0)=g(x_0)=0$, in a certain neighborhood of the point
$x=x_0,\varepsilon=0$ the functions $F(x,\varepsilon)$ and
$G(x,\varepsilon)$ are smooth functions.

{\bf Theorem 1.} {\it If the inequalities
$$\begin{aligned}
&f^{\prime\prime}(x_0)g^\prime(x_0)-g^{\prime\prime}(x_0)f^\prime(x_0)<0,\\
&g'(x_0)>0,\,\,\,F(x_\varepsilon,\varepsilon)>0,\end{aligned}\eqno(11)
$$
where $x_\varepsilon$ is a zero of the function $G(x,\varepsilon)$ in
the neighborhood of the point  $x=x_0$, are satisfied,
then there exists a number $\varepsilon_0>0$ such that for
all $\varepsilon\in(0,\varepsilon_0)$ system (10) has a cycle.}

{\bf Proof.} We introduce the notation $z(t)=x(t)-x_0$ and put
$$\begin{aligned}
&f(x)=f_1(x-x_0)+f_2(x-x_0)^2+O((x-x_0)^3),\\
&g(x)=g_1(x-x_0)+g_2(x-x_0)^2+O((x-x_0)^3).\end{aligned}
$$
Without loss of generality, it can be assumed that $g_1=1$.

Consider the case when $\varepsilon=0$ and $z(0)$ is a small number.

Here we shall use smoothnis of functions $F$ and $G$ and shall
follow the first Lyapunov method on finite time interval [Lefschetz,
1957; Cesari, 1959].

The first approximation    to the solution $z(t)$ of system (10)
is the function
$$
z_1(t)=z(0)\cos t.
$$
The second approximation $z_2(t)$ to the solution $z(t)$ is given by
the equation
$$
\ddot z_2+z_2=u(t),\eqno(12)
$$
where
$$
u(t)=z(0)^2(f_1\cos t\sin t-g_2(\cos t)^2).
$$
By Lemma 1 we obtain that $z_2(t)$ has the form
$$\begin{aligned}
&z_2(t)=z(0)\cos t+z(0)^2\left(-\frac{f_1}{3}(\sin t)^3\cos t+\right.\\
&+\frac{g_2}{3}(1-(\cos t)^3)\cos t+\frac{f_1}{3}(1-(\cos t)^3)\sin
t-\\
&-\left.g_2(\sin t-\frac{1}{3}(\sin t)^3)\sin
t\right).\end{aligned}\eqno(13)
$$
This implies that for $\dot z(0)=0$ the relation holds:
$$\begin{aligned}
&\dot z_2(t)=-z(0)\sin t+z(0)^2\left(\frac{f_1}{3}(\sin t)^4\right.-\\
&-\frac{g_2}{3}(1-(\cos t)^3)\sin t+\frac{f_1}{3}(1-(\cos t)^3)\cos
t-\\
&-g_2\left.(\sin t-\frac{1}{3}(\sin t)^3)\cos t\right).\end{aligned}
$$
It follows that the time  $T>0$, of the second crossing, of the
straight line $\{y=0\}$ by the solution $x(t), y(t)$ of system (10)
with the initial data $x(0),\,\,y(0)=0$  satisfies the relation
$$
T=2\pi+O((x(0)-x_0)^2).\eqno(14)
$$

Consider now the function
$$
V(x,y)=\left(y+\intl_{x_0}^xf(z)dz\right)^2+2\intl_{x_0}^xg(z)dz.
$$

 For the derivative of the function $V$ along
the solution of system (10)  the following relation
$$\begin{aligned}
\dot V(x,&y)=-2g(x)\intl_{x_0}^xf(z)dz=\\
&=-f_1(x-x_0)^3-\left(\frac{2}{3}f_2+f_1g_2\right)(x-x_0)^4+\\
&+O((x-x_0)^5)\end{aligned}
$$
is valid.
 Then, taking into account relations (13), (14), we obtain
$$\begin{aligned}
&V(x(T),y(T))-V(x(0),0)=\\
&=-\intl_0^T\left(f_1z_2(t)^3+\left(\frac{2}{3}f_2+f_1g_2\right)z_2(t)^4\right)dt+\\
&+O(z(0)^5)=-z(0)^4\intl_0^{2\pi}\left(\frac{2}{3}f_2+f_1g_2\right)(\cos
t)^4+\\
&+3f_1(\cos t)^2\left(-\frac{f_1}{3}(\sin t)^3\cos
t+\frac{g_2}{3}(1-\right.\\
&-\left.(\cos t)^3\right)\cos t+\frac{f_1}{3}(1-(\cos t)^3)\sin t-\\
&-g_2(\sin t-\left.\frac{1}{3}(\sin t)^3)\sin t\right)dt+\\
&+O(z(0)^5)=-\frac{(f_2-f_1g_2)\pi}{2}z(0)^4+O(z(0)^5).\end{aligned}
$$
  In this case by the theorem on a continuous dependence of solutions
of system (10) on a parameter, for sufficiently small parameter
$\varepsilon>0$     in comparison with the number $|z(0)|$ we have
the inequality
$$
V(x(T),y(T))>V(x(0),0).
$$
 On the other hand, for small $\varepsilon>0$
the equilibrium $x=x_\varepsilon, y=0$ of system (10) is a stable
focus. These two facts yield  that in a certain small neighborhood
of the point $x=x_\varepsilon, y=0$  there exists a cycle of system
(10) (Fig. 1).
\begin{figure}[!ht]
\begin{center}
 \includegraphics[scale=0.8]{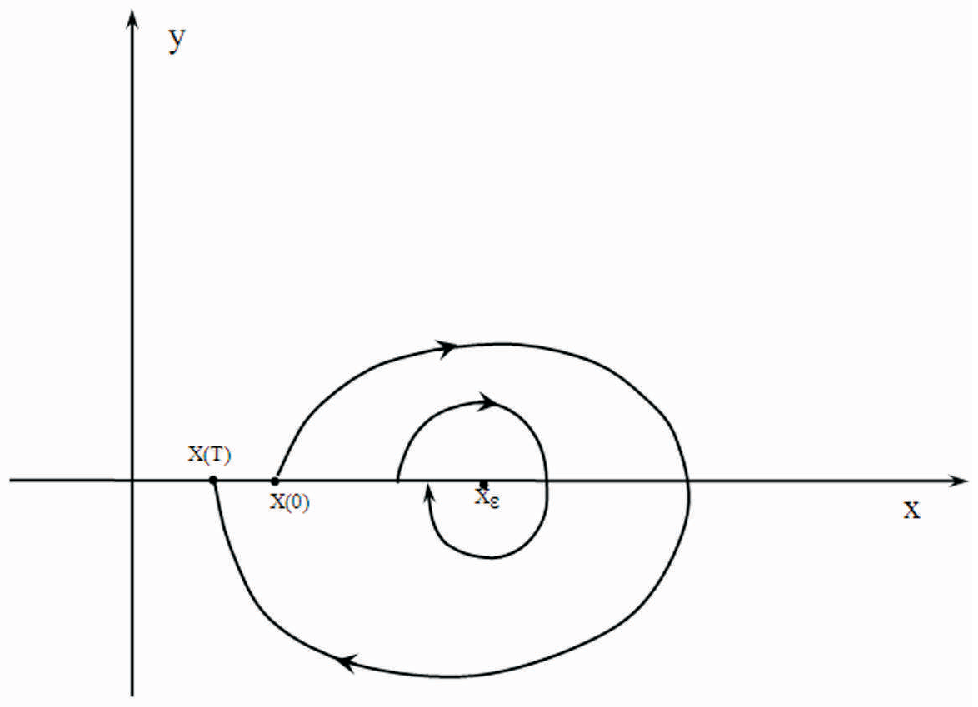}
\caption{}
\end{center}
\end{figure}
Theorem is proved. $\Box$

{\bf Example 1.} Consider  the equation
$$
\ddot x+(\varepsilon+f_1x+f_2x^2)\dot x+x+g_2x^2+g_3x^3=0,
$$
 where $f_1,f_2,g_2,g_3$ are arbitrary numbers, $f_2-f_1g_2<0$.

Theorem 1 implies that for the sufficiently  small positive
$\varepsilon$ this equation has at least one periodic
solution.

Consider system (1) and assume that
$$\begin{aligned}
&b_1=\beta_1=1,\,\alpha_1+\beta_2=\varepsilon,\\
&A=-b_2+2a_1c_2-a_1,\\
&B=-b_2-\beta_2+2\alpha_1c_2-2a_1,\\
&C=-a_2-2\alpha_2+\alpha_1b_2+a_1\beta_2+\alpha_1\beta_2-c_2\alpha_1^2,\\
&D=-\alpha_2+\alpha_1\beta_2.\end{aligned}
$$

If $D>0$ and
$$
AD-BC+BD(1+c_2)<0
$$
then by theorem 1 we obtain that for sufficiently small
$\varepsilon>0$ the system (1) has a cycle in some neighborhood of
the point $x=y=0$.

Suppose that in system (1)  $a_1=a_2=c_1=0$, $b_1=\beta_1=c_2=1$,
$b_2=-1,\,\alpha_1=1\!/3-\varepsilon,\,\alpha_2=\beta_2=-1\!/3$. We
have
$$\begin{aligned}
&F(x,\varepsilon)=(x^2+2(1-\varepsilon)x+\varepsilon)|1+x|^{-3},\\
&G(x,\varepsilon)=\left(\frac{1}{9}(x^2+2x)+\varepsilon x(x^2+2x+\frac{1}{3}-\varepsilon x)\right)\\
&(1+x)^{-5},\\
&f'(0)=f'(-2)=2,\,\,g'(0)=g'(-2)=\frac{2}{9}\\
&f''(0)=f''(-2)=-10,\,\,g''(0)=g''(-2)=-2.\end{aligned}
$$

In the neighborhood of $x=0$, $x_\varepsilon=0$
and in the neighborhood of $x=-2$ the relation
$x_\varepsilon=-2-3\varepsilon+o(\varepsilon)$ holds.

For $x_\varepsilon=0$, we have
$$
F(x_\varepsilon,\varepsilon)=\varepsilon
$$
and for $x_\varepsilon=-2-3\varepsilon+o(\varepsilon)$,
$$
F(x_\varepsilon,\varepsilon)=\frac{11\varepsilon+o(\varepsilon)}{|1+x_\varepsilon|^3}.
$$

It should also be noted that in the considered case   system (10)
with $\varepsilon=0$    under the change
$$
x=-z-2,\,\,\,y=-\vartheta
$$
preserves its form and the equilibrium $x=-2,\,y=0$    goes into
the equilibrium $z=\vartheta=0$.

Thus, we prove the following

{\bf Theorem 2.} {\it For sufficiently small $\varepsilon<0$, the
system
$$\begin{aligned}
&\dot x=xy+\left(\frac{1}{3}-\varepsilon\right)x+y\\
&\dot y=-xy+y^2-\frac{1}{3}x-\frac{1}{3}y\end{aligned}
$$
has    no less than two cycles. Moreover, each of two equilibria is
surrounded by   no less than one cycle.}

\bigskip

\begin{center}{\bf 4.   The generalization of the Lienard theorem}\end{center}

\bigskip

In this section we give the generalization [Leonov, 2006] of a
classical result of Lienard concerning  the existence of periodic
solution for the equation [Lefschetz, 1957]
$$
\ddot x+f(x)\dot x+g(x)=0\eqno(15)
$$

Consider a system
$$\begin{aligned}
&\dot x=y\\
&\dot y=-f(x)y-g(x),\end{aligned},\eqno(16)
$$
which is equivalent to equation (15).

Suppose, the function $f(x)$ and $g(x)$ are continuous on the
interval $(a,+\infty)$ and for the certain numbers $a<\nu_1\le x_0\le
\nu_2$ the following conditions

 $$\begin{aligned} {\rm 1)}\quad &\lim\limits_{x\to a}g(x)=-\infty,\\
 &\lim\limits_{x\to +\infty}g(x)=+\infty,\\
 &\lim\limits_{x\to a}\int\limits_{x_0}^xg(z)dz=
 \lim\limits_{x\to
 +\infty}\int\limits_{x_0}^xg(z)dz=+\infty,\end{aligned}\eqno(17)
 $$
 $$\begin{aligned}
{\rm 2)}\quad & f(x)>0,\quad\forall\,
x\in(a,\nu_1)\bigcup(\nu_2,+\infty),\\
&\int\limits_{\nu_2}^{\nu_1}f(x)dx\le 0\end{aligned}\eqno(18)
$$
are satisfied.

{\bf Theorem  3.} {\it If conditions 1) and 2) are valid, then in
the phase space
$$
\{x\in(a,+\infty),\,y\in R^1\}\eqno(19)
$$
there exists a piecewise-smooth transversal closed curve, which
intersects the straight line $\{y=0\}$ at the certain points
$a<\mu_1<\nu_1$ and $\mu_2>\nu_2$ (Fig. 2). If, in addition, in
space (19) system (16) has only one unstable by Lyapunov focal
equilibrium, then system (16) has a cycle.}

{\bf Proof.} Consider a pair of the numbers $\mu_1\in(a,\nu_1)$ and
$\mu_2\in(\nu_2,+\infty)$ such that $\mu_1$ is sufficiently close to
$a$, $\mu_2$ is sufficiently large, and
$$
\int\limits_{\mu_1}^{\mu_2}g(x)dx=0.\eqno(20)
$$

Without loss of generality, we put
$$
g(x)<0,\quad \forall\, x\in[\mu_1,\nu_1],
$$
$$
g(x)>0,\quad \forall\, x\in[\nu_2,\mu_2].
$$

Consider the following functions
$$\begin{aligned}
& V_1(x,y)=y^2+2\int\limits_{x_0}^xg(z)dz,\\
& V_2(x,y)=\left(y+\int\limits_{\nu_1}^xf(z)dz\right)^2+2\int\limits_{x_0}^xg(z)dz,\\
& V_3(x,y)=\left(y+\int\limits_{\nu_2}^xf(z)dz\right)^2+2\int\limits_{x_0}^xg(z)dz,\\
&V_4(x,y)=V_2(x,y)-\varepsilon(x-\nu_1)\\
&V_5(x,y)=V_3(x,y)-\varepsilon(x-\nu_2)\\
&V_6(x,y)=V_3(x,y)+\varepsilon(x-\nu_2)\\
&V_7(x,y)=V_2(x,y)+\varepsilon(x-\nu_1).\end{aligned}
$$
Here $\varepsilon$ is a certain sufficiently small number.

 Consider now the sets $\Omega_j$ (see Fig. 2)
$$\begin{aligned}
&\Omega_1=\{x\in[\mu_1,\nu_1], y\ge 0, V_1(x,y)\le V_1(\mu_1,0)\},\\
&\Omega_2=\{x\in[\nu_1,x_0], y\ge 0, V_4(x,y)\le V_2(\nu_1,y_1)\},\\
&\Omega_3=\{x\in[x_0,\nu_2], y\ge 0, V_5(x,y)\le V_3(\nu_2,y_2)\},\\
&\Omega_4=\{x\in[\nu_2,\mu_2], y\ge 0, V_3(x,y)\le V_3(\mu_2,0)\},\\
&\Omega_5=\{x\in[\nu_2,\mu_2], y\le 0, V_1(x,y)\le V_1(\mu_2,0)\},\\
&\Omega_6=\{x\in[x_0,\nu_2], y\le 0, V_6(x,y)\le V_3(\nu_2,y_3)\},\\
&\Omega_7=\{x\in[\nu_1,x_0], y\le 0, V_7(x,y)\le V_2(\nu_1,y_4)\},\\
&\Omega_8=\{x\in[\mu_1,\nu_1], y\le 0, V_2(x,y)\le V_2(\mu_1,0)\}.
\end{aligned}
$$

\begin{figure}[!ht]
\begin{center}
 \includegraphics[scale=0.3]{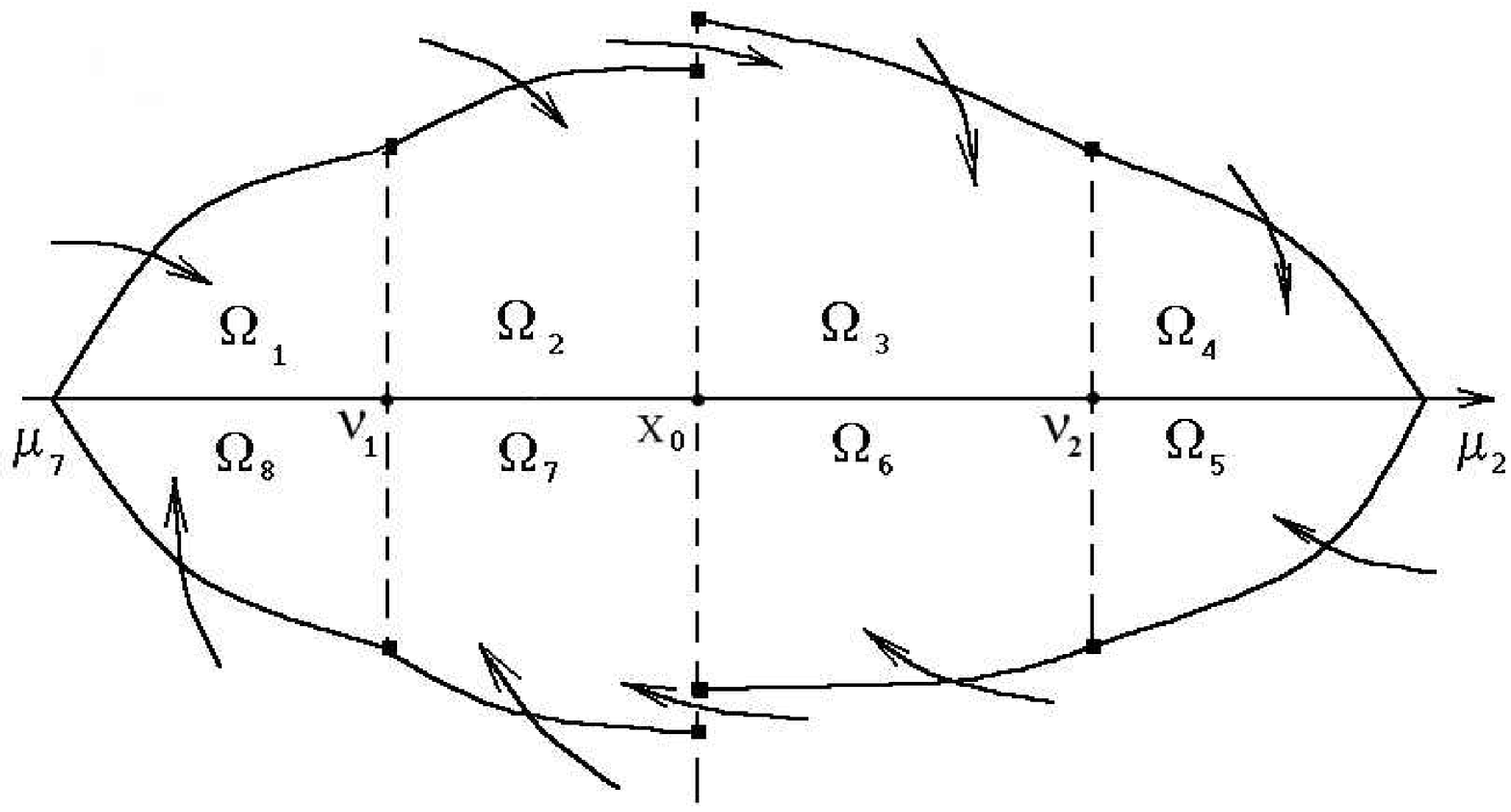}
\caption{}
\end{center}
\end{figure}

Here $y_1>0, y_2>0, y_3<0, y_4<0$ are  solutions of the square equations
$$\begin{aligned}
&V_1(\nu_1,y_1)=V_1(\mu_1,0),\\
&V_3(\nu_2,y_2)=V_3(\mu_2,0),\\
&V_1(\nu_2,y_3)=V_1(\mu_2,0),\\
&V_2(\nu_1,y_4)=V_2(\mu_1,0)\end{aligned}
$$

For the derivatives of the functions $V_j(x,y)$ along the solutions
of system (24)  we have the following relations
$$\begin{aligned}
&\dot V_1=-2f(x)y^2,\quad \dot
V_2=-2g(x)\int\limits_{\nu_1}^xf(z)dz,\\ &\dot
V_3=-2g(x)\int\limits_{\nu_2}^xf(z)dz,\\
&\dot V_4=-2g(x)\int\limits_{\nu_1}^xf(z)dz-\varepsilon y,\\
&\dot V_5=-2g(x)\int\limits_{\nu_2}^xf(z)dz-\varepsilon y,\\
&\dot V_6=-2g(x)\int\limits_{\nu_2}^xf(z)dz+\varepsilon y,\\ &\dot
V_7=-2g(x)\int\limits_{\nu_1}^xf(z)dz+\varepsilon y.\end{aligned}
$$

Therefore under the above assumptions    concerning $\mu_1$ and $\mu_2$
and for $y\ne 0, x\ne\nu_j$ we have the inequalities

$$\begin{aligned}
&\dot V_1<0\quad\mbox{on}\quad \Omega_1\cup\Omega_5,\quad \dot
V_4<0\quad\mbox{on}\quad \Omega_2,\\
&\dot V_5<0\quad\mbox{on}\quad \Omega_3,\quad \dot
V_3<0\quad\mbox{on}\quad \Omega_4,\\
&\dot V_6<0\quad\mbox{on}\quad \Omega_6,\quad \dot
V_7<0\quad\mbox{on}\quad \Omega_7,\\
&\dot V_2<0\quad\mbox{on}\quad \Omega_8.\end{aligned}
$$

Note that by conditions  2), (20)  for sufficiently small
$\varepsilon$ we have inequalities $y_5<y_6$ and $y_7>y_8$, where
$y_5$ is a positive solution of the equation
$$
V_4(x_0,y_5)=V_2(\nu_1,y_1),
$$
$y_6$ is a positive solution of the equation
$$
V_5(x_0,y_6)=V_3(\nu_2,y_2),
$$
$y_7$ is a negative solution of the equation
$$
V_6(x_0,y_7)=V_3(\nu_2,y_3),
$$
$y_8$ is a negative solution of the equation
$$
V_7(x_0,y_8)=V_2(\nu_1,y_4).
$$
(Fig. 2).

Thus, we construct a transversal closed curve, which proves the
first assertion of theorem. The instability of focal equilibrium
gives the existence of cycle. $\Box$

 Let us clarify the fact that  the relations
$\dot V_4<0, \dot V_5<0, \dot V_6<0,\dot V_7<0$ are satisfied on the sets
$\Omega_2, \Omega_3,\Omega_6,\Omega_7$, respectively.

Hold fixed the arbitrary $\varepsilon>0$, we choose $\mu_1$ and
$\mu_2$ so much closer to $a$ and to $+\infty$, respectively, that
the minimal values of $|y|$ on the intersection of the closed curve
(Fig. 2) and the band $\{x\in[\nu_1,\nu_2]\}$ are more than
$$
\frac{1}{\varepsilon}\max\limits_{x\in[\nu_1,\nu_2]}2
\left|g(x)\int\limits_{\nu_1}^xf(z)dz\right|
$$
and
$$
\frac{1}{\varepsilon}\max\limits_{x\in[\nu_1,\nu_2]}2
\left|g(x)\int\limits_{\nu_2}^xf(z)dz\right|.
$$
This implies the required inequalities $\dot V_j<0$. Thus, we have
here $\mu_j=\mu_j(\varepsilon)$  and
$$
\lim\limits_{\varepsilon\to 0}\mu_1(\varepsilon)=a,\quad
\lim\limits_{\varepsilon\to 0}\mu_2(\varepsilon)=+\infty.
$$

Represent now conditions 1) and 2) in terms of quadratic system (1).

We have
$$
a=-\frac{\beta_1}{b_1},\quad b_1\ne 0.
$$
Without loss of generality, we put $b_1>0$. Conditions 1) and 2)
are satisfied if
$$
 0<2c_2<b_1,\quad\beta_1>0,
$$
$$
\frac{a_1\beta_1}{b_1}>\alpha_1
$$
$$
\frac{a_1(2c_2-b_1)}{b_1}>b_2\eqno(21)
$$
$$
\frac{a_1(b_1b_2-a_1c_2)}{b_1^2}>a_2.
$$

Besides, from conditions (21) it follows a positive invariance of
the half-plane $\{x\ge a\}$ for quadratic system (1).

Note that here $c_1=0$ and the parameters $\alpha_2$ and $\beta_2$
do not enter into conditions (21).

 A set of closed piecewise-smooth transversal curves can be constructed
similarly. Therefore we have the following

{\bf Theorem  4.} {\it Let conditions (21) be valid. Then any
solution, of system (1) with initial data such that $x(0)>a$, tends
as $t\to+\infty$ to a bounded
 attractor, placed in the half-plane $\{x>a\}$.}

Theorem  4 permits us to localize the search of the limit cycles of
quadratic system (1). We remark that if in the half-plane $\{x>a\}$
we have a unique unstable by Lyapunov focal equilibrium of system
(1) and for it  conditions (21) are valid, then system (1) has a
periodic solution, placed in the half-plane $\{x>a\}$.

Besides, under these assumptions we have $f(x)>0$ for $x<a$.
Then, using the  Lyapunov functions
$$V(x,y)=y^2+\int\limits_{2a}^xg(z)dz
$$
 we can prove that a solution of system (1) with  the initial data
$x(0)<a$ either tends as $t\to+\infty$ to the     equilibrium,
either to infinity, either leaves in a definite time the half-plane
$\{x<a\}$.

Using the above and Theorem  4 we

{\bf Theorem  5.} {\it Suppose, conditions (21) are valid and in the
half-plane $\{x>a\}$ system (1) has a unique unstable by Lyapunov
focal equilibrium. Then any solution of system (1)  with initial
data such that $x(0)<a$ either tends as $t\to+\infty$ to the
equilibrium, either to infinity, either to a bounded attractor,
situated in the half-plane $\{x>a\}$. Any solution of system (1)
with initial data such that $x(0)\ge a$ tends as $t\to+\infty$ to a
bounded attractor, situated in the half-plane $\{x>a\}$. This
attractor has at least one cycle.}

{\bf Example 2.}  We put $b_1=\beta_1=1,\,\alpha_1=0,\,
a_2=b_2=-1,\,
c_1=0,\,c_2=1\!/4,\,\alpha_1=-1,\,\beta_2=2,\,\alpha_2=-1000$. Here
conditions (21) are fulfilled, $g(x)\ne 0$ for $x\ne 0$ and $x>-1$,
equilibrium $x=y=0$ is unstable focus. Therefore system (1) has at
least one cycle by theorem 5.

The conditions of Theorem  5    effectively select in the space of
parameters of system (1) a set, of positive Lebesque measure, in
which the cycles exist. What the obtaining of the results of such
kind is of present interest is remarked in [Arnol'd, 2005].

V.A.Arnold (see [Arnol'd, 2005]) writes: "To estimate the number of
limit cycles of quadratic vector fields on a plane, A.N.Kolmogorov
distributed, as a mathematical practice, few hundreds of such fields
(with the randomly chosen coefficients of quadratic expressions)
   among few hundreds students of the mechanics and mathematics faculty of
Moscow State University.

 Each student  had to find the number of limit cycles of his field.

The result of this experiment was    perfectly    sudden:  all
fields  had    none of the limit cycles!

For a small change of  coefficients of field, a limit cycle is
preserved. Therefore in the space of coefficients, the systems with
one, two, three (and even, as was proved later, four) limit cycles
make up open sets and if the coefficients of polynomials are drawn
at random, then the    probability of hitting in this sets are
positive.

The fact that  this event did not occur    suggests that the
above-mentioned probabilities are, evidently, small."

But the set of existence of cycles with selected by theorem 5 is not
small.

\bigskip

\noindent{\bf References}

\medskip

\noindent Arnol'd V.I. [1976] {\it Ordinary Differential Equations}

(Nauka, Moscow) (in Russian).

\noindent Arnol'd V.I. [2005] {\it Experimental mathematics} \break
\indent (Fazis, Moscow) (in Russian).

\noindent Blows T.R., Perko L.M. [1994] ``Bifurcation of limit

cycles from centers and separatrix cycles of

planar analytic systems", {\it SIAM Review}, {\bf 36},\break \indent
N 3,341--376.

\noindent Cesari L. [1959] {\it Asymptotic Behavior and Stabi\-lity

Problems in Ordinary Differential  Equations}

 (Springer, Berlin).

\noindent Chen Lan Sun, Wang Mind Shu. [1979] `` The

relative position and the number of limit \break \indent cycles of a
quadratic differential systems", \break \indent
 {\it Acta Math.
Sinica}, {\bf 22}, 751--758.

\noindent Hilbert D. [1902] ``Mathematical problems", {\it Bull.

Amer. Math. Soc.}, {\bf 8}, 437--479.

\noindent Ilyashenko Yu. [ 2002] ``Centennial history of \break
\indent Hilbert's 16th problem Bulletin of the AMS",\break
\indent{\bf 39}, N 3. 301--354.

\noindent Lefschetz S. [1957] {\it Differential Equations:
Geo-

metric Theory} (Interscience Publishers, New

York).

\noindent Leonov G.A. [1997] ``Two-Dimensional Quadra-\break\indent
tic Systems as a Lienard Equation", {\it Differen-\break\indent
tial
Equations and Dynamical Systems}, {\bf 5},

N 3/4, 289--297.

\noindent Leonov G.A. [1998] ``The problem of estimation of

the number of cycles of two-dimensional quadra-\break\indent
tic
systems from nonlinear mechanics point of

view", {\it Ukr. Math. J.}, {\bf 50}, N 1, 48--57.

\noindent Leonov G.A. [2006] ``Family of transversal curves

for two-dimensional dynamical systems",

{\it Vestnik St.Petersburg University}, {\bf 1}, N 4, 48--78.

\noindent Lloyd N.G. [1988] ``Limit cycles of polynomial
sys-\break\indent tems -- some recent developments", {\it In book:

New Direction in Dynamical Systems.

 Cambridge University Press}, 192--234.

\noindent Shi Song Ling. [1980] ``A concrete example of the

existence of four limit cycles for plane quadratic

systems", {\it Sci. Sinica}, {\bf 23},  153--158.

\end{document}